\documentclass[graybox]{svmult}

\usepackage{mathptmx}       
\usepackage{helvet}         
\usepackage{courier}        
\usepackage{type1cm}        
                            
\usepackage{makeidx}        
\usepackage{graphicx}       
                            
\usepackage{multicol}       
\usepackage[bottom]{footmisc}

\usepackage{amsfonts}
\usepackage{enumerate}
\usepackage{amsmath}
\usepackage{amssymb}
\usepackage{amsxtra}
\usepackage{amscd}
\usepackage{colonequals}
\usepackage{url}

\makeindex 

\begin{document}

\title*{Elementary deformations and the hyperK\"ahler-quaternionic
K\"ahler correspondence}
\titlerunning{Elementary deformations and the hK/qK correspondence}
\author{Oscar Macia and Andrew Swann}
\institute{Oscar Macia
\at
Departamento de Geometria y Topologia,
Facultad de Ciencias Matematicas,
Universidad de Valencia,
C. Dr. Moliner, 50,
Burjassot (46100) Valencia,
Spain
\email{oscar.macia@uv.es} 
\and Andrew Swann
\at
Department of Mathematics,
Aarhus University,
Ny Munkegade 118,
Bldg 1530,
DK-8000 Aarhus C,
Denmark
\email{swann@imf.au.dk}}  

\maketitle

\abstract{The hyperK\"ahler-quaternionic K\"ahler correspondence
constructs quaternionic K\"ahler metrics from hyperK\"ahler metrics
with a rotating circle symmetry.  We discuss how this may be
interpreted as a combination of the twist construction with the
concept of elementary deformation, surveying results of our forthcoming
paper.  We outline how this leads to a uniqueness statement for the
above correspondence and indicate how basic examples of c-map
constructions may be realised in this context.}

\section{Introduction}
\label{sec:introduction}

The twist construction was introduced in \cite{Swann:T,Swann:twist} as
a geometric construction that reproduces T-duality arguments in the
physicists literature for geometries with torsion.  It has proved
successful in constructing compact simply connected examples of a
number of classes of non-K\"ahler geometries.  However, elsewhere in
the physics literature string theory dualities are used to construct
metrics of special holonomy.  In particular, the c-map construction of
Cecotti et al.~\cite{Cecotti-FG:II} produces quaternionic K\"ahler
metrics from projective special K\"ahler manifolds.  An intermediate
stage in this construction is a passage from hyperK\"ahler manifolds
of a given to dimension to quaternionic K\"ahler manifolds of the same
dimension.

HyperK\"ahler and quaternionic K\"ahler metrics are two of the
infinite families of geometries in the holonomy classification of
Berger~\cite{Berger:hol,Besse:Einstein}.  They are both Einstein
geometries and there are many open questions about their structure and
classification.  In 2008, Haydys \cite{Haydys:hKS1} showed how to each
quaternionic K\"ahler manifold with circle action one may associate
hyperK\"ahler manifolds with a symmetry that fixes only one of the
complex structures.  He also provided a description of how to invert
that construction.  Later Hitchin \cite{Hitchin:hK-qK} gave a twistor
interpretation of this construction along the lines of
\cite{Joyce:hypercomplex}, and
\cite{Alexandrov-DP:qK-hK,Alekseevsky-CM:conification} provided
expressions in arbitrary signature.  The metric constructions here all
have the flavour of making a conformal change, but with two different
factors along and orthogonal to directions determined by a symmetry.

The purpose of this note is to describe the results of
\cite{Macia-S:c-map}, where we show that the twist construction can be
used to interpret this so-called hyperK\"ahler-quaternionic K\"ahler
correspondence at to prove that there is only one degree of freedom
this construction.  We then indicate how the computational framework
of the twist construction may be applied to understand some of the
basic examples of the c-map.

\section{Twist constructions}
\label{sec:twist-constructions}

The twist construction~\cite{Swann:T,Swann:twist} associates to a
manifold with a circle action a new space of the same dimension with a
distinguished vector field.  

Suppose \( M \) is manifold of dimension \( n \).  Let \( X \) be a
vector field on \( M \) that generates a circle action.  A twist~\( W
\) of \( M \) is specified as a quotient \( W = P/\langle X'\rangle \)
of a principal \( S^1 \)-bundle~\( P \to M \) by a lift \( X' \) of~\(
X \).  It thus fits in to a double fibration
\begin{equation*}
  \begin{CD}
    M @<{\pi_M}<< P @>{\pi_W}>> W.
  \end{CD}
\end{equation*}
If \( H^2(M,\mathbb Z) \) has no torsion, the bundle \( P \) is
specified by the curvature form~\( F \) of a connection one-form \(
\theta \in \Omega^1(P) \), given by \( \pi_M^*F = d\theta \).  We let
\( \mathcal H = \ker\theta \) be the corresponding horizontal
distribution on~\( P \).  Lifts \( X' \) of \( X \) that preserve \(
\theta \) and the principal vector field \( Y \) are given by
\begin{equation*}
  X' = X^\theta + (\pi_M^*a)Y,
\end{equation*}
where \( X^\theta \in \mathcal H \) is the horizontal lift of~\( X \)
with respect to~\( \theta \) and \( a \in C^\infty(M) \) is a
Hamiltonian function satisfying
\begin{equation}
  \label{eq:da}
  da = - X {\lrcorner} F.
\end{equation}
This requires that \( F \) is preserved by~\( X \).  The \emph{twist}
\( W \colonequals P/\langle X'\rangle \) then
admits a circle action generated by~\( (\pi_W)_*Y \).

This essentials of this set-up are specified by the \emph{twist data}
\( (M,X,F,a) \) with \( X \in \mathfrak X(M) \) generating a circle
action, \( F \in \Omega_{\mathbb Z}^2(M)^X \) an \( X \)-invariant
closed two-form with integral periods and \( a \)
satisfying~\eqref{eq:da}.

Provided \( a \) is non-zero, invariant tensors on \( M \) may be
transferred to \( W \) as follows.  Note that at \( p\in P \), the
projections \( \pi_M \) and \( \pi_W \) induce isomorphisms \(
T_{\pi_M(p)}M \cong \mathcal H_p \cong T_{\pi_W(p)}W \).  Thus given
\( p\in \pi_M^{-1}(q) \), a tensor \( \alpha_q \) at \( q \in M \)
induces a tensor \( (\alpha_W)_{q'} \) at \( q' = \pi_W(p) \in W \).
The tensor \( \alpha_W \) is well-defined if \( \alpha \) is preserved
by \( X \).  We then say that \( \alpha \) and \( \alpha_W \) are \(
\mathcal H \)-related and write
\begin{equation*}
  \alpha \sim_{\mathcal H} \alpha_W.
\end{equation*}

The two most important computational facts for \( \mathcal H
\)-related tensors are:

\begin{property}
  for \( \alpha \in \Omega^p(M)^X \) an invariant \( p \)-form, the
  exterior differential of \( \alpha_W \) is given by
  \begin{equation}
    \label{eq:dW}
    d\alpha_W \sim_{\mathcal H} d_W\alpha \colonequals d\alpha -
    \frac1aF\wedge X {\lrcorner} \alpha.
  \end{equation}
\end{property}

\begin{property}
  for an invariant complex structure \( I \) on \( M \) that is \(
  \mathcal H \)-related to an almost complex structure \( I_W \) on
  \( W \), we have
  \begin{equation*}
    \text{\( I_W \) is integrable if and only if \( F \) is of type \(
    (1,1) \) for~\( I \).}
  \end{equation*}
\end{property}

\noindent
Recall that \( F \in \Omega^2(M) \) is of type \( (1,1) \) if \(
F(IA,IB) = F(A,B) \) for all \( A,B \in TM \).

These facts show that geometric properties of the twist are determined
by the twist data.

\begin{example}
  A basic example of the twist construction is provided by \( M =
  \mathbb CP(n) \times T^2 \).  This is a K\"ahler manifold as a
  product.  Suppose \( X \) generates one of the circle factors of \(
  T^2 = S^1 \times S^1 \).  Taking \( F \) to be the Fubini-Study
  two-form on \( \mathbb CP(n) \), we have \( X {\lrcorner} F = 0 \),
  so can take \( a \equiv 1 \).  Then \( P = S^{2n+1} \times T^2 \)
  and the twist is \( W = S^{2n+1} \times S^1 \).  As \( F \)~is
  type~\( (1,1) \) we have that \( W \) is a complex manifold.
  However \( W \)~is compact and \( b_2(W) = 0 \), so \( W \) can not
  be K\"ahler.
\end{example}

\section{Elementary deformations of hyperK\"ahler metrics}
\label{sec:elem-deform}

As formula~\eqref{eq:dW} indicates, the twist of a closed differential
form is rarely closed.  In a given geometric situation it is therefore
interesting to adjust the geometric data before performing a twist.

We wish to work with hyperK\"ahler manifolds.  These are
(pseudo-)Riemannian manifolds \( (M,g) \) with almost complex
structures \( I \), \( J \) and~\( K \) such that
\begin{enumerate}
\item \( IJ = K = -JI \),
\item \( g \) is Hermitian with respect to \( I \), \( J \) and \( K \),
\item the two-forms \( \omega_I = g(I\cdot,\cdot) \), \( \omega_J \)
  and \( \omega_K \) are closed:
  \begin{equation*}
    d\omega_I = 0 = d\omega_J = d\omega_K.
  \end{equation*}
\end{enumerate}
By Hitchin~\cite{Hitchin:Riemann-surface} the last condition implies
that \( I \), \( J \) and \( K \) are integrable. The restricted
holonomy is then a subgroup of \( \mathit{Sp}(n) \), where \( \dim M =
4n \), and the metric is Ricci-flat.  The triples \( (g,I,\omega_I)
\), etc., are then K\"ahler structures on~\( M \).

Let \( X \) be a symmetry of a hyperK\"ahler structure \( (M,g,I,J,K)
\), but which we mean that \( X \) is an isometry that preserves the
linear span \( \langle I,J,K \rangle \) of \(
I,J,K\in\operatorname{End}(TM) \).  The vector field \( X \) induces
four one-forms on \( M \) given by
\begin{equation*}
  \alpha_0 = g(X,\cdot),\quad \alpha_I = I\alpha_0 =
  -\alpha(I\cdot),\quad \alpha_J = J\alpha_0,\quad \alpha_K = K\alpha_0.
\end{equation*}
We then define
\begin{equation*}
  g_\alpha = \alpha_0^2 + \alpha_I^2 + \alpha_J^2 + \alpha_K^2.
\end{equation*}
When \( X \) is not null, \( g_\alpha \) is positive semi-definite
proportional to the restriction of~\( g \) to \( \mathbb HX = \langle
X,IX,JX,KX \rangle \).

\begin{definition}
  An \emph{elementary deformation} of a hyperK\"ahler metric \( g \)
  with respect to a symmetry \( X \) is a metric of the form
  \begin{equation*}
    g^N = fg + hg_\alpha
  \end{equation*}
  with \( f \) and \( h \) smooth functions on~\( M \).
\end{definition}

\noindent This is thus more general than a conformal change of \( g \).

As \( I \), \( J \) and \( K \) are parallel, we have that \( X \)
acts as a linear transformation on \( \mathbb R^3 = \langle I,J,K
\rangle \).  It preserves the algebraic relations, so acts as an
element of~\( \mathfrak{so}(3) \).  As \( \mathfrak{so}(3) \) has rank
one, it follows that the action is either trivial or conjugate a
circle action fixing \( I \) and mapping \( J \) to~\( K \).  By
relabelling the complex structures and normalising~\( X \) we may thus
assume in this latter case that
\begin{equation}
  \label{eq:XJ}
  L_XI = 0\quad\text{and}\quad L_XJ = K.
\end{equation}
An isometry \( X \) satisfying \eqref{eq:XJ} will be called
\emph{rotating}. 

For a rotating symmetry, we have \( d\alpha_I = 0 \), \( d\alpha_J =
\omega_K \) and \( d\alpha_0 = G - \omega_I \), where \( G \in
\Omega^2(M) \) is a two-form that is of type \( (1,1) \) for \( I \),
\( J \) and~\( K \).  As \( \alpha_I \) is closed, we may pass to a
cover of \( M \) and write \( \alpha_I = d\mu \) for a smooth map \(
\mu\colon M \to \mathbb R \).  The function  \( \mu \) is a K\"ahler
moment map for \( X \) with respect to \( (g,\omega_I) \).

\section{The hyperK\"ahler-quaternionic K\"ahler correspondence}
\label{sec:hyperk-quat-kahl}

Suppose \( (M,g,I,J,K) \) is hyperK\"ahler with rotating symmetry~\( X
\) with K\"ahler moment map~\( \mu \).  Then \( X \) does not preserve
\( \omega_J \) or \( \omega_K \), but the four-form
\begin{equation}
  \label{eq:Omega}
  \Omega = \omega_I^2 + \omega_J^2 + \omega_K^2
\end{equation}
is invariant and closed.

If \( W \) is manifold of dimension at least~\( 8 \) with a four-form
\( \Omega^W \) pointwise linearly equivalent to~\eqref{eq:Omega}, then
\( W \) has an almost quaternion-Hermitian structure \( (g_W,\mathcal
G) \), where \( \mathcal G\subset \operatorname{End}(TW) \) is a
three-dimensional subbundle with a local basis \( (I_W,J_W,K_W) \) of
almost complex structures for which \( g_W \) is Hermitian and with \(
I_WJ_W=K_W= -J_WI_W \).  Such a structure is \emph{quaternionic
K\"ahler} if \( \Omega^W \) is parallel with respect to the
Levi-Civita connection of~\( g_W \).  It follows that \( g_W \) is
Einstein \cite{Salamon:Invent,Besse:Einstein}.  If \( \dim W\geqslant
12 \), then to obtain quaternionic K\"ahler it is sufficient that \(
d\Omega^W = 0 \) \cite{Swann:symplectiques}.  For \( \dim W = 8 \),
one can work with the local two-forms \( \omega^W =
(\omega^W_I,\omega^W_J,\omega^W_K) \) and quarternionic K\"ahler is
then equivalent to the existence of a local connection form \( A \in
\Omega^1(\mathfrak{so}(3)) \) such that \( d\omega^W = A\wedge\omega^W
\).

The behaviour of the exterior derivative under the twist is given by
\eqref{eq:dW}, so from the above remarks we may determine whether a
twist is quaternionic K\"ahler by working on~\( M \).

\begin{theopargself}
  \begin{theorem}[\cite{Macia-S:c-map}]
    Let \( (M,g,I,J,K) \) by a hyperK\"ahler manifold with non-null
    rotating symmetry~\( X \) and K\"ahler moment map~\( \mu \).  If
    \( \dim M \geqslant 8 \) then, up to homothety, the only twists of
    elementary deformations~\( g^N = fg + hg_\alpha \) of \( g \) that
    are quaternionic K\"ahler have
    \begin{equation}
      \label{eq:gN}
      g^N = \frac1{(\mu-c)^2}g_\alpha - \frac1{\mu-c}g
    \end{equation}
    for some constant~\( c \).  Furthermore the corresponding twist
    data is given by
    \begin{equation*}
      F = kG = k(d\alpha_0 + \omega_I),\quad a = k(g(X,X)-\mu+c),
    \end{equation*}
    for some constant \( k \).
  \end{theorem}
\end{theopargself}

The method of proof is first to impose the quaternionic K\"ahler
condition on as arbitrary twist of~\( \Omega^N \), the four-form
associated to \( g^N \) via \( I \), \( J \) and \( K \), and to
decompose these equations in type components relative to \( \mathbb HX
\) and its orthogonal complement.  From this one deduces that \( f \)
and \( h \) are functions of \( \mu \) and that \( h = f' \).  Then we
consider the equation \( da = - X {\lrcorner} F\) and determine the twist
function~\( a \).  Finally, we investigate the condition \( dF = 0 \),
which provides an ordinary differential equation for~\( f \).

From the Theorem, it follows that the constructions provided in
\cite{Haydys:hKS1,Hitchin:hK-qK,Alekseevsky-CM:conification} of
quaternionic K\"ahler metrics from hyperK\"ahler metrics with rotating
circle symmetry agree.

\begin{example}
  \label{ex:flat}
  We consider \( \mathbb H^{p,q} = \mathbb R^{4p,4q} \), \( n=p+q \),
  with its flat hyperK\"ahler metric
  \begin{equation*}
    g = \sum_{i=1}^n \varepsilon_i(dx_i^2 + dy_i^2 + du_i^2 + dv_i^2)
  \end{equation*}
  with \( \varepsilon_i = +1 \), for \( i\leqslant p \), and \(
  \varepsilon_i = -1 \), for \( i>p \), %
  and K\"ahler two-forms
  \begin{gather*}
    \omega_I = \sum_{i=1}^n \varepsilon_i(dx_i \wedge dy_i - du_i
    \wedge dv_i),\quad \omega_J = \sum_{i=1}^n \varepsilon_i(du_i
    \wedge dx_i + dv_i
    \wedge dy_i),\\
    \omega_K = \sum_{i=1}^n \varepsilon_i(du_i \wedge dy_i - dv_i
    \wedge dx_i).
  \end{gather*}
  If \( X \) is a rotating circle symmetry then it is an element of \(
  \mathfrak{sp}(p,q) + \mathfrak{u}(1) \), but lies in a maximal
  compact subgroup, so is conjugate to
  \begin{equation*}
    X = \sum_{i=1}^n
    \bigl(\frac12-\lambda_i\bigr)\bigl(y_i\frac\partial{\partial x_i}
    -x_i\frac\partial{\partial y_i}\bigr) -  \bigl(\frac12 +
    \lambda_i\bigr)\bigl(v_i\frac\partial{\partial u_i} -
    u_i\frac\partial{\partial v_i}\bigr)
  \end{equation*}
  for some \( \lambda_1,\dots,\lambda_n \in \mathbb R \).  For \( X \)
  to be non-null, we must have \( \sum_{i=1}^n\varepsilon_i\lambda_i^2
  \ne 0 \).  This vector field has \( d\alpha_0 = d(g(X,\cdot)) = G -
  \omega_I \) with
  \begin{equation*}
    G = 2\sum_{i=1}^n \varepsilon_i\lambda_i(dx_i\wedge dy_i +
    du_i\wedge dv_i) 
  \end{equation*}
  so \( G = d\beta \), where \( \beta = \sum_{i=1}^n
  \varepsilon_i\lambda_i(-y_idx_i+x_idy_i - v_idu_i + u_idv_i) \).

  The twisting form \( F \) is equal to a multiple of \( G=d\beta \),
  so is exact and the twist bundle is trivial \( P = \mathbb H^n
  \times S^1 \).  Let us take \( F = G \).  The connection one-form
  may be written as \( \theta = \beta + d\tau \), where
  \( \partial/\partial\tau \) generates the principal \( S^1
  \)-action.  The horizontal lift \( X^\theta \) of \( X \) to \( P \)
  is then
  \begin{equation*}
    X^\theta = X - \beta(X)\frac\partial{\partial\tau}.
  \end{equation*}
  Direct calculation shows that \( d(\beta(X)) = - X {\lrcorner} F \),
  so the twist function is \( a = \beta(X) + c \) and the twist is the
  quotient of \( P \) by \( X' = X + c\frac\partial{\partial\tau} \).
  Thus the twist is
  \begin{equation*}
    W = (\mathbb H \times S^1)/ \bigl\langle X +
    c\frac\partial{\partial\tau} \bigr\rangle.
  \end{equation*}
  This will be an orbifold if \( \lambda_i \) and \( c \) are
  integers.  It is smooth when they are pairwise co-prime.

  The theorem says that \( W \) is equipped with a quaternionic
  K\"ahler metric \( \mathcal H \)-related to~\( g^N \) in
  equation~\eqref{eq:gN}, whenever this is non-degenerate.  The
  function \( \mu \) is given by
  \begin{equation*}
    \mu = \frac12 \sum_{i=1}^n \bigl(\frac12-\lambda_i\bigr)(x_i^2 +
    y_i^2) + \bigl(\frac12+\lambda_i\bigr)(u_i^2+v_i^2).
  \end{equation*}
  The metric \( g^N \) has two contributions to its signature.  On the
  quaternionic span~\( \mathbb HX \) of~\( X \), the sign is that of
  \( (\lVert X \rVert^2-\mu+c)\lVert X \rVert^2/(\mu-c)^2 \),
  orthogonal to \( \mathbb HX \) the original metric is multiplied by
  \( -\lVert X \rVert^2/(\mu-c) \).  Thus up to overall sign \( g^N \)
  has quaternionic signature that is either \( (p+1,q-1) \), \( (p,q)
  \) or \( (p-1,q+1) \).  It is degenerate on the sets \( (\lVert X
  \rVert^2 = 0) \), i.e., where \( X \) is null, and on \( (\lVert X
  \rVert^2-\mu+c = 0) \), which is the set where the twist function \(
  a \) vanishes.  The metric may also blow-up on \( (\mu = c) \).
\end{example}

\section{Application to the c-map}
\label{sec:applications-c-map}

The c-map is a construction introduced by Cecotti et
al.~\cite{Cecotti-FG:II}.  It starts with a so-called projective
special K\"ahler manifold of dimension \( 2n \) and produces a
quaternionic K\"ahler manifold of dimension \( 4n+4 \).  Explicit
local expressions for the resulting metrics where provided by Ferrara
and Sabharwal~\cite{Ferrara-S:q}.  Recently Alekseevsky et
al.~\cite{Alexandrov-DP:qK-hK,Alekseevsky-CDM:qK-special} have shown
that the hyperK\"ahler-quaterionic K\"ahler correspondence reproduces
the quaternionic K\"ahler metrics of the c-map.  In particular, this
means that one may obtain all the known examples homogeneous (positive
definite) quaternionic K\"ahler of negative scalar curvature, and
their work is also beginning to produce new examples of complete
quaternionic K\"ahler metrics.

Given the wide generality of the twist construction, is useful to
understand how such homogeneous examples may arise.  To be concrete
let us consider the real hyperbolic space \( \mathbb RH(2) \) as a
solvable Lie group with K\"ahler metric of constant curvature.  This
has a global basis~\( \{a,b\} \) of one forms, such that \( da = 0 \)
and \( db = - \lambda\alpha\wedge b \), for some constant \( \lambda
\) depending on the scalar curvature.  For this to be a projective
special K\"ahler manifold, we need to consider a certain cone metric
and show that it admits a flat symplectic connection of special
K\"ahler type, as described by Freed~\cite{Freed:special}.  

Let \( C_0 \) be a circle bundle over \( \mathbb RH(2) \) with
connection one-form \( \varphi \) whose curvature is \( 2a\wedge b
\).  Pulling \( a \) and \( b \) back to \( C = \mathbb R_{>0} \times
C_0 \), the cone geometry is
\begin{equation*}
  g_C = t^2(a^2+b^2-\varphi^2) - dt^2,\quad \omega_C = t^2 a\wedge b -
  t\varphi \wedge dt,
\end{equation*}
a K\"ahler metric of signature \( (2,2) \).  It has a symmetry \( X \)
generated by the principal action on~\( C_0 \).

Locally, one can show that this admits a special K\"ahler connection
if and only if \( \lambda ^2 \) is \( 4 \) or \( 4/3 \).  In case \(
\lambda^2 = 4 \), the special connection agrees with the Levi-Civita
connection of~\( g_c \).  In both cases, using the cotangent
trivialisation \( (\hat a,\hat b,\hat\varphi,\hat\psi) =
(ta,tb,t\varphi,dt) \), one may construct a hyperK\"aher metric of
signature \( (4,4) \) on \( H = T^*C \) of the form \( g_H = \hat a^2
+ \hat b^2 - \hat \varphi^2 - \hat\psi^2 + \hat A^2 + \hat B^2 - \hat
\Phi^2 - \hat \Psi^2 \).  Indeed the flat connection gives \( TH = V^*
\oplus V \), with \( V \cong TM \).  This is the rigid c-map, see
Freed~\cite{Freed:special}.  The K\"ahler forms \( \omega_J \) and \(
\omega_K \) are just the real and imaginary parts of the standard
complex symplectic two-form on \( H = T^*C \).

Horizontally lifting the symmetry of~\( X \) of \( C \) to~\( H=T^*C
\) using the flat connection, one obtains a rotating symmetry \(
\widetilde X \) of the hyperK\"ahler structure.  Note that the
symmetry \( X \) does not preserve the flat connection, and it rotates
the quadruple \( \tilde\delta = (\hat A,\hat B,\hat \Phi, \hat \Psi)
\).  The twist data for this lifted action is given by the curvature
form
\begin{equation*}
  F = - \hat a\wedge\hat b + \hat \varphi\wedge\hat\psi -\hat
  A\wedge\hat B + \hat \Phi\wedge\hat \Psi
\end{equation*}
and twist function~\( -t^2/2 + c \).  The curvature form is exact, and
so we may proceed much as in Example~\ref{ex:flat}.  

In particular, we have a coordinate~\( \tau \) on~\( S^1 \) in \( P =
H \times S^1 \).  With \( c=0 \) the twist is then diffeomorphic to \(
(H/\langle \widetilde X \rangle) \times S^1 \).  We may use \( \tau \)
to define a new quadruple \( \delta = \tilde\delta \exp(\mathbf i\tau)
\), where \( \mathbf i = \operatorname{diag}(\mathbf i_2,\mathbf i_2)
\) with \( \mathbf i_2 = \left(
  \begin{smallmatrix}
    0&-1\\ 1&0
  \end{smallmatrix}
\right) \).  Now using \eqref{eq:dW} one may show that the structure
functions of the coframe \( \mathcal H \)-related to \( (\hat a,\hat
b, \hat\varphi, \hat\psi, \delta_1, \delta_2, \delta_3, \delta_4) \)
are constants, so these define a dual basis for a Lie algebra.  The
metric \( g^N \) is seen to be positive definite, complete and has
constant coefficients in this coframe, so the resulting quaternionic
K\"ahler metric on~\( W \) is complete.  It follows that the universal
cover of \( W \) is a Lie group~\( G \) and that the metric on \( W \)
pulls back to a left-invariant metric.  We have \( W = G/\mathbb Z \)
and knowing the structure constants we may identify \( G \) as the
solvable Lie groups that act transitively on the non-compact symmetric
spaces \( \operatorname{Gr}_2(\mathbb C^{2,2}) \) for \( \lambda^2=4
\) or \( G_2^*/\mathit{SO}(4) \) for \( \lambda^2=4/3 \).  This
provides a global verification of the main example of Ferrara and
Sabharwal~\cite{Ferrara-S:q}.

\begin{acknowledgement}
This work is partially supported by the Danish Council for
Independent Research, Natural Sciences, and by the  Spanish Agency for 
Scientific and Technical Research (DGICT) and FEDER project
MTM2010-15444.
 
\end{acknowledgement}

\end{document}